%% file: preprint.tex
\documentclass{article}
\usepackage{arxiv}

\usepackage[utf8]{inputenc} 
\usepackage[T1]{fontenc}    
\usepackage{microtype}      
\usepackage{doi}

\usepackage{amsmath,amsfonts}
\usepackage{algorithmic}
\usepackage{algorithm}
\usepackage{array}
\usepackage{textcomp}
\usepackage{stfloats}
\usepackage{url}
\usepackage{verbatim}
\usepackage{graphicx}
\usepackage{cite}

\usepackage{bbm}
\usepackage{hyperref}
\usepackage{siunitx}
\DeclareSIUnit{\euro}{\text{€}}
\DeclareSIUnit{\year}{\text{a}}
\DeclareSIUnit{\voltamperereactif}{\text{var}}

\usepackage{mathtools}
\usepackage{tikz,pgfplots}
\usepgfplotslibrary{fillbetween,groupplots,patchplots}
\usetikzlibrary{patterns,math,spy}
\pgfplotsset{samples=600,compat=newest,every axis/.append style={font=\footnotesize}}
\usepackage{tabularx}
\usepackage{booktabs}
\usepackage{makecell}
\usepackage{multirow}
\usepackage[capitalise]{cleveref}

\usepackage{enumitem}

\input{shortcuts}

\title{Improving Operational Feasibility in Large-Scale Power System Planning}
\date{}

\usepackage{authblk}

\newbox{\orcid}\sbox{\orcid}{\includegraphics[scale=0.06]{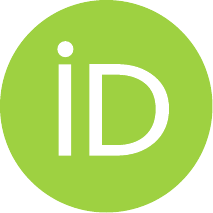}} 
\author[1,2]{%
	\href{https://orcid.org/0009-0000-6286-147X}{\usebox{\orcid}\hspace{1mm}Gereon~Recht\thanks{\texttt{gereon.recht@dlr.de}}}%
}
\author[1,2]{%
	\href{https://orcid.org/0000-0002-0080-7862}{\usebox{\orcid}\hspace{1mm}Oussama~Alaya}%
}
\author[2]{%
	\href{https://orcid.org/0009-0003-4504-5165}{\usebox{\orcid}\hspace{1mm}Benedikt~Jahn}%
}
\author[1]{%
	\href{https://orcid.org/0000-0002-9720-0337}{\usebox{\orcid}\hspace{1mm}Karl-Kiên~Cao}%
}
\author[2]{%
	\href{https://orcid.org/0000-0002-0208-4100}{\usebox{\orcid}\hspace{1mm}Hendrik~Lens}%
}
\affil[1]{German Aerospace Center (DLR), Institute of Networked Energy Systems, Curiestr. 4, 70563 Stuttgart, Germany}
\affil[2]{Institute of Energy Process Engineering and Dynamics in Energy Systems (IED), University of Stuttgart, Germany}

\hypersetup{
pdftitle={Improving Operational Feasibility in Large-Scale Power System Planning},
pdfauthor={Gereon~Recht, Oussama~Alaya, Benedikt~Jahn, Karl-Kiên~Cao, Hendrik~Lens},
}

\begin{document}

\maketitle

\begin{abstract}
    Large-scale power system planning mostly uses linearized, active power only approximations of the power flow equations, ignores many operational constraints, and tests the operational feasibility of the resulting systems only under strongly simplifying assumptions.
    We propose an approach to obtain solutions to large instances of the alternating current capacity expansion problem via redispatch and reinforcement of an initial solution.
    The problem formulation considers simultaneous expansion of generators, reactive compensation devices, storage systems, and transmission.
    Furthermore, it includes operational constraints via startup procedures and capability curves of power sources and simplified stability limits via constraints on voltage angle differences and voltage magnitudes.
    To obtain initial solutions, we test several established and partly modified power flow approximations and integrate them into an approach for iterative transmission expansion planning, thereby obtaining convex formulations.
    We demonstrate the approach on large problem instances covering the islands of Great Britain and Ireland at the transmission level, for which we extend the open data source to model reactive power.
    We find that including transmission losses to determine the initial solution is most decisive, as the amount of redispatch and reinforcements necessary to obtain an alternating current feasible solution is reduced, whereas incorporating reactive power constraints did not lead to further improvements.
    Our approach ensures an alternating current feasible system under weak assumptions, thus guaranteeing steady-state voltage stability and allowing subsequent dynamic grid simulations, which is instrumental for planning stable future inverter-dominated power systems.
\end{abstract}

\section*{Nomenclature}
\begin{description}[leftmargin=!,labelwidth=\widthof{$\underline{z}_l = r_l + j x_l$}]
    \item[$\Tset$] set of snapshots $\Tset = \{1,\dots,T\}$
    \item[$\Nset$] set of electrical buses
    \item[$\Lset$] set of branches
    \item[$\LsetRev$] set of reversed branches 
    \item[$\LsetAC, \LsetDC \subseteq \Lset$] sets of high voltage alternating current and high voltage direct current branches
    \item[$\mathcal{C}$] cycle basis of network $(\Nset, \LsetAC)$
    \item[$\Sset$] set of power sources
    \item[$\SsetSg, \SsetIbr \subseteq \Sset$] sets of synchronous generators and inverter-based resources
    \item[$\SsetSto \subseteq \Sset$] set of storage systems
    \item[$\Sset_n \subseteq \Sset$] set of power sources attached to bus $n$
    \item[$\mathcal{K}_s^{\text{up}}, \mathcal{K}_s^{\text{lo}}$] sets of tuples defining lines for approximating PQ-capability curves 
    \item[$\mathcal{P}_{l,0}$] set of points at which the tangents of a linear relaxation are centred 
    \item[$\delta_{t}$] duration of snapshot in hours
    \item[$d_{n,t}^p, d_{n,t}^q$] active and reactive load
    \item[$\underline{z}_l = r_l + \jim x_l$] branch impedance
    \item[$\underline{y}_l = g_l + \jim b_l$] branch admittance
    \item[$b_l^{\text{sh}}$] branch shunt susceptance
    \item[${f}_l^{\text{max}}$] thermal limit of branch
    \item[$\theta_l^{\text{max}}$] max.\ voltage angle difference across branch
    \item[$\eta_l$] losses of high voltage direct current branch
    \item[$c_s, c_l$] capital costs for power sources and transmission lines
    \item[$o_s, o_s^{\text{su}}$] marginal and startup costs of power source
    \item[$a_{s,t}$] per unit availability of power source
    \item[$p_{s}^{\text{min}}, p_{s}^{\text{max}}$] active power injection limits
    \item[$q_{s}^{\text{min}}, q_{s}^{\text{max}}$] reactive power injection limits
    \item[$p_{s,t}^{\text{in},\text{max}}$] max.\ inflow of energy into storage system
    \item[$e_{s}^{\text{max}}$] max.\ state of charge of storage system
    \item[$\eta_s, \eta_s^c$] discharging and charging efficiencies of storage system
    \item[$K, C$] incidence and cycle incidence matrices 
    \item[$\theta_{n,t}$, $v_{n,t}$] voltage angle and magnitude
    \item[$\theta_{l,t}$] voltage angle difference $\theta_{n,t} - \theta_{m,t}$ over branch $l=(n,m)$
    \item[$u_s$, $u_l$] number of power source units and transmission line circuits
    \item[$p_{l,t}$, $q_{l,t}$] active and reactive power flows at bus $n$ of branch $l = (n,m)$
    \item[$p_{l,t}^{\text{loss}}$, $q_{l,t}^{\text{dem}}$] approximations of losses and reactive power demand of branch
    \item[$\cosrelax_{l,t}$] variable representing relaxation of $\cos \theta_{l,t}$
    \item[$\beta_{s,t}$] number of online units of power source
    \item[$p_{s,t}$, $q_{s,t}$] active and reactive power injection 
    \item[$p^c_{s,t}$] charging of storage system 
    \item[$p_{s,t}^{\text{in}}$] controllable inflow into storage system 
    \item[$e_{s,t}$] state of charge of storage system 
\end{description}

\section{Introduction}
To date, large instances of capacity expansion problem formulations that include the equations describing steady-state alternating current (AC) power flow and unit commitment are computationally intractable.
This is unlikely to change fundamentally due to the NP-hardness of both integer programming and the AC power flow problem \cite{karp_reducibility_1972,BIENSTOCK2019494}.
Hence, it is common practice to allow continuous expansion decisions, ignore unit commitment, and use the linearized, active power only DC power flow approximation, which may include an approximation of transmission losses \cite{NEUMANN2022118859}.
If the AC-feasibility of the resulting system is investigated, strong assumptions such as infinite capacities for voltage control and the existence of a slack generator with infinite capacity are used \cite{NEUMANN2022118859}.
Moreover, no guarantees on the AC-feasibility of the resulting systems can be given.
This is not sufficient for prospective studies that aim to consider steady-state voltage-stability or conduct a posteriori dynamic stability analyses, which require AC operating points that take into account the operational constraints of all power sources \cite{kundur}.
Power flow approximations offering improvements over the DC approximation have been used for power system planning in numerous other studies \cite{WOGRIN2020115925,DAVE2021106683,nagarajan2017resilienttransmissiongriddesign,6308747,7038299}.
These novel approximations may improve the deviation of the objective value from the optimal AC solution, but generally do not produce AC-feasible solutions \cite{VENZKE2020106480}.
Common approaches to determine AC-feasible solutions include load shedding \cite{nagarajan2017resilienttransmissiongriddesign,us_tep} and constraint tightening \cite{7038299,DAVE2021106683}. 
Another possibility is reinforcing a solution of the approximated problem, which we call \emph{initial solution} in the following, such that it becomes AC-feasible \cite{6308747}.
With the reasoning that the initial solution is expected to be somewhat \emph{near} an AC-feasible solution, a practical motivation for this approach is the existence of heuristic solvers like IPOPT \cite{wachter_implementation_2006}, which aim to determine local optima and are known to reliably find solutions to comparatively small benchmark instances of the AC optimal power flow (\textsc{ac-opf}) problem \cite{VENZKE2020106480}.
However, the approximation error of the heuristic solution is not bounded and to date, a polynomial time approximation scheme exists only for a special case of \textsc{ac-opf} \cite{8892494}.

With the exception of \cite{NEUMANN2022118859}, where only active power was considered, all mentioned publications only contain applications on small benchmark problem instances.
To date, there are no studies on the applicability of advanced power flow approximations to large-scale instances with simultaneous planning of generator, storage, and transmission expansion.
Moreover, to the best of the authors' knowledge, there are no studies that determine an AC-feasible solution for such large problem instances.
We hence provide a formulation of an AC capacity expansion problem which includes constraints on operational and stability limits and on reactive power.
By integrating several power flow approximations, we obtain different convex formulations of the problem.
In a reinforcement approach, we are able to determine an AC-feasible system based on an initial solution by redispatching or installing additional power sources and reactive power compensation devices.
We evaluate the suitability of each power flow approximation to provide an initial solution to this approach.
The methodology is demonstrated on highly resolved problem instances with tens of millions of variables and constraints covering Great Britain and Ireland at the transmission level.
Furthermore, we compile a set of assumptions for considering reactive power in studies based on open data.

\section{Methodology}\noindent
We first define the intractable \textsc{alternating current capacity expansion problem} (\textsc{ac-cep}) and four simplified variants of this problem, each using a different power flow approximation.
A summary of the constraints defining each variant is given in \cref{tab:model_definitions}.
We investigate differences regarding modelling of transmission losses in the power flow approximations and present an approach to reinforce a solution for AC-feasibility.
A potential issue of using voltage angle difference constraints within an established successive convex programming approach for considering the expansion-related change of transmission line parameters is discussed.
\subsection{The capacity expansion problem}\noindent
The \textsc{cep} aims to minimize the total system costs which comprises capital costs for expansion of generators, storage systems and branches, as well as operational costs.
The resulting system has to be feasible for a sequence of snapshots $t \in \Tset = \{1,\dots,T\}$ of varying load and availability of generation.
Snapshots can represent time spans of different durations $\delta_t$.
We assume a given network $(\Nset, \Lset)$ of electrical buses $\Nset$ and branches $\Lset$.
Power sources $\Sset$ and branches $\Lset$ are counted in units $u_s$ and circuits $u_l$, which have associated annualised capital cost of $c_s$ and $c_l$, respectively.
We allow continuous expansion within $u_s \in [u_s^{\text{min}}, u_s^{\text{max}}]$, $u_l \in [u_l^{\text{min}}, u_l^{\text{max}}]$.
The operational costs $o_s$ define the costs per unit of electrical energy injected into the bus.
Within the set of power sources $\Sset$, we differentiate between synchronous generators (SGs) $\SsetSg$ and inverter-based resources (IBRs) $\SsetIbr$.
The variable $\beta_{s,t}$ represents the number of online units, which can be increased via the startup variable $\beta_{s,t}^{\text{su}}$ at a cost of $o_s^{\text{su}}$.
We define the objective
\begin{equation} \label{eqn:obj}
  \min \sum_{s \in\Sset} \left( c_s u_s + \sum_{t \in \Tset} \left( o_s p_{s,t} \delta_{t} + o_s^{\text{su}} \beta_{s,t}^{\text{su}} \right)  \right) + \sum_{l \in \Lset} c_l u_l,
\end{equation}
which corresponds to minimizing the annualised system costs.
We model the behaviour of power sources as in \cite{11305559} via
\begin{align}
  & 0 \leq \beta_{s,t} \leq u_s & \forall s \in \SsetSg, t \in \Tset \label{eqn:online_limit}, \\
  & \beta_{s,t} \leq \beta_{s,t-1} + \beta_{s,t}^{\text{su}} & \forall s \in \SsetSg, t \in \Tset \setminus \{1\} \label{eqn:online_consistency}, \\
  & \beta_{s,1} \leq \beta_{s,T} + \beta_{s,1}^{\text{su}} & \forall s \in \SsetSg \label{eqn:online_boundary}.
\end{align}
While $\beta_{s,t}$ is not necessary for IBRs, as they are not governed by the same operational limits as SGs, we will use it for both to simplify notation.
Active and reactive power injections are constrained via
\begin{align}
  & p^{\text{min}}_s \beta_{s,t} \leq p_{s,t} \leq a_{s,t} p^{\text{max}}_s \beta_{s,t} & \forall s \in \Sset, t \in \Tset \label{eqn:ps_injections_online}, \\
  & q^{\text{min}}_s \beta_{s,t} \leq q_{s,t} \leq q^{\text{max}}_s \beta_{s,t} & \forall s \in \Sset, t \in \Tset \label{eqn:qs_injection_online}.
\end{align}
Storage systems are additionally associated with the variables $p_{s,t}^c$, $e_{s,t}$ for charging and the state of charge, respectively.
Their behaviour is constrained by
\begin{align}
  & 0 \leq p^c_{s,t} \leq p^{\text{max}}_s \beta_{s,t} & \forall s \in \SsetSto, t
  \in \Tset \label{eqn:storage_unit_charging}, \\
  & p_{s,t} + p^c_{s,t} \leq p_s^{\text{max}} \beta_{s,t} & \forall s \in \SsetSto, t
  \in \Tset \label{eqn:storage_unit_complementarity_relaxation}, \\
  & 0 \leq e_{s,t} \leq e^{\text{max}}_s u_s & \forall s \in \SsetSto, t
  \in \Tset \label{eqn:state_of_charge_limits}, \\
  & e_{s,t} = e_{s,t-1} + \delta_t \left(\eta^c_s p^c_{s,t} - \frac{p_{s,t}}{\eta_s} + p_{s,t}^{\text{in}} \right) & \forall s \in \SsetSto, t \in \Tset \setminus \{1\}, \label{eqn:state_of_charge} \\
  & e_{s,1} = e_{s,T} & \forall s \in \SsetSto \label{eqn:state_of_charge_cyclic},
\end{align}
as in \cite{11305559}, where we include a relaxation of the complementarity constraint in \labelcref{eqn:storage_unit_complementarity_relaxation} to limit simultaneous charging and discharging of the same storage unit. 
We approximate the PQ-capabilities of each power source as a convex polytope via a set of constraints
\begin{align}
  p_{s,t} & \leq \tau_{s,k}^{\text{up}} q_{s,t} + \upsilon_{s,k}^{\text{up}} p^{\text{max}}_s \beta_{s,t} & \forall s \in \Sset, t \in \Tset, (\tau_{s,k}^{\text{up}}, \upsilon_{s,k}^{\text{up}}) \in \mathcal{K}_s^{\text{up}},  \label{eqn:pq_upper} \\
  p_{s,t} & \geq \tau_{s,k}^{\text{lo}} q_{s,t} + \upsilon_{s,k}^{\text{lo}} p^{\text{max}}_s \beta_{s,t} & \forall s \in \Sset, t \in \Tset, (\tau_{s,k}^{\text{lo}}, \upsilon_{s,k}^{\text{lo}}) \in \mathcal{K}_s^{\text{lo}}, \label{eqn:pq_lower}
\end{align}
where the parameters $\tau_{s,k}$, $\upsilon_{s,k}$ represent the slope and per unit intercept of lines acting as upper and lower bounds for $p_{s,t}$.

High voltage alternating current (HVAC) transmission lines are represented by their $\Pi$-equivalent circuit, which models a line $l \in \LsetAC$ using a series admittance $\underline{y}_l = g_l + \jim b_l$ and a shunt susceptance $b_l^{\text{sh}}$ that is distributed equally to both buses.
To represent power flows, we introduce variables and constraints for both endpoints of each HVAC branch.
To that end, we define the set of reversed branches $\LsetRev$ such that $(n, m) \in \Lset \Leftrightarrow (m, n) \in \LsetRev$.
The index $l^\prime$ then refers to the reversed counterpart of branch $l$.
We now model the active and reactive power flow at bus $n \in \Nset$ of branch $l = (n,m) \in \LsetAC \cup \LsetACRev$ using the variables $p_{l,t}$ and $q_{l,t}$, respectively.
Given the bus voltage angles $\theta_{n,t}, \theta_{m,t} \in [-\pi/2, \pi/2]$ and the per unit voltage magnitudes $v_{n,t}, v_{m,t} \in [0.9, 1.1]$, the flows are given by
\begin{align}
    p_{l,t} =& g_l v_{n,t}^2 - v_{n,t} v_{m,t} g_l \cos\theta_{l,t} - v_{n,t} v_{m,t} b_l \sin \theta_{l,t} & \forall l \in \LsetAC \cup \LsetACRev, t \in \Tset, \label{eqn:ac_flows_p} \\
    q_{l,t} =& - (b_l + \frac{b_l^{\text{sh}}}{2}) v_{n,t}^2 - v_{n,t} v_{m,t} g_l \sin \theta_{l,t} + v_{n,t} v_{m,t} b_l \cos\theta_{l,t} & \forall l \in \LsetAC \cup \LsetACRev, t \in \Tset, \label{eqn:ac_flows_q}
\end{align}
where $\theta_{l,t} = \theta_{n,t} - \theta_{m,t}$ represents the voltage angle difference across the branch.
Each branch has a thermal limit
\begin{align}\label{eqn:thermal_limit}  
  &p_{l,t}^2 + q_{l,t}^2 \leq (a_l f^{\text{max}}_l u_l)^2 &\forall l \in \LsetAC \cup \LsetACRev, t \in \Tset
\end{align}
where the parameter $a_l \in [0, 1]$ represents a limit for the maximum line loading.
Voltage angle differences are limited via
\begin{align}\label{eqn:voltage_angle_difference}
  &-\theta^{\text{max}}_l \leq \theta_{l,t} \leq \theta^{\text{max}}_l & \forall l \in \LsetAC, t \in \Tset
\end{align}
to ensure normal operation within a specified maximum voltage angle difference.
Eq.\ \labelcref{eqn:voltage_angle_difference} also represents a technical requirement for power flow approximations in the following, which can include approximations or relaxations that are only defined within a certain range of $\theta_{l,t}$.
Also, tighter voltage angle difference limits often improve the approximation error \cite{lpac,qc_relaxation}.

We model the flow over high voltage direct current (HVDC) branches with constant losses $\eta_l$ depending linearly on the branch length as in \cite{NEUMANN2022118859}.
To apply the losses, we model the bidirectional flow using two variables $\overrightarrow{p}_l, \overleftarrow{p}_l \geq 0$ representing the flows along and against the orientation of the branch such that the power flows at the sending and receiving bus are given by
\begin{align}
  p_{l,t} &= \overrightarrow{p}_{l,t} - (1 - \eta_l)  \overleftarrow{p}_{l,t} &\forall l \in \LsetDC, t \in \Tset, \label{eqn:flow_hvdc_along}   \\
  p_{l^\prime,t} &= \overleftarrow{p}_{l,t} - (1 - \eta_l) \overrightarrow{p}_{l,t} &\forall l' \in \LsetDCRev, t \in \Tset. \label{eqn:flow_hvdc_against}
\end{align}
The flow is limited by
\begin{align}\label{eqn:flow_hvdc_limit}
  &\vert p_{l,t} \vert \leq p^{\text{max}}_l u_l & \forall l  \in \LsetDC \cup \LsetDCRev, t \in \Tset.
\end{align}

We differentiate between HVDC branches connected via line-commutating converters and voltage-source converters.
For line-commutating converters, we neglect the reactive power demand of converter stations by assuming that it is compensated locally and thus contained in the capital costs \cite{kundur}.
Voltage-source converters can provide reactive power compensation, which we model via power sources attached to the sending and receiving endpoint of a HVDC connection whose injection is limited by the capacity of the associated HVDC branch.
Naturally, there is no flow of reactive power across HVDC branches, i.\,e.\ it holds that $q_{l,t} = 0$ for all $l \in \LsetDC \cup \LsetDCRev$.

Finally, we enforce Kirchhoff's current law via the nodal balance equations
\begin{align}
  &\sum_{s \in\Sset_n} p_{s,t} - \sum_{s \in \mathcal{S}_n^{\text{sto}}} p_{s,t}^c - d_{n,t}^p =\sum_{\substack{l=(n,m) \\ \in \Lset \cup \LsetRev}} p_{l,t} & \forall n \in \Nset, t \in \Tset, \label{eqn:ac_nodal_balance_p} \\
  &\sum_{s \in\Sset_n} q_{s,t} + d_{n,t}^q = \sum_{\substack{l=(n,m) \\ \in \Lset \cup \LsetRev}} q_{l,t} & \forall n \in \Nset, t \in \Tset, \label{eqn:ac_nodal_balance_q}
\end{align}
where we neglect bus shunts and assume constant power loads, i.\,e.\ loads that are voltage-independent.

Since we have relaxed the unit commitment and investment decisions, the remaining sources of nonconvexity in \textsc{ac-cep} are given by the power flow constraints in \labelcref{eqn:ac_flows_p,eqn:ac_flows_q} and the dependence of the line parameters on the number of circuits \cite{8916411}, rendering \textsc{ac-cep} NP-hard \cite{BIENSTOCK2019494}.
In the following sections, we discuss how we resolve these nonconvexities via approximations of the power flow equations and a successive convex programming approach for transmission expansion planning.

\subsubsection{DC approximation}\noindent
The DC approximation linearises the power flow equations through a well-known set of assumptions \cite{192975}.
All reactive power constraints (\labelcref{eqn:ac_flows_q,eqn:ac_nodal_balance_q,eqn:pq_lower,eqn:pq_upper,eqn:qs_injection_online,eqn:ac_flows_p}) are ignored and \labelcref{eqn:ac_flows_p,eqn:ac_nodal_balance_p,eqn:thermal_limit,eqn:voltage_angle_difference,eqn:ac_flows_p} are replaced by
\begin{align}
  & \sum_{l \in \LsetAC} C_{lc} x_l p_{l,t} = 0 & \forall c \in \Cset, t \in \Tset, \label{eqn:kvl} \\
  & \vert p_{l,t} \vert \leq a_l f_l^{\text{max}} u_l & \forall l \in \LsetAC, t \in \Tset, \label{eqn:thermal_limit_dc} \\
  & \vert p_{l,t} \vert \leq \frac{\theta^{\text{max}}_l}{x_l} & \forall l \in \LsetAC, t \in \Tset \label{eqn:voltage_angle_difference_dc_approx}, \\
  &\sum_{s \in\Sset_n} p_{s,t} - \sum_{s \in \mathcal{S}_n^{\text{sto}}} p_{s,t}^c - d_{n,t}^p = \sum_{l \in \LsetAC} K_{nl} p_{l,t} + \sum_{\substack{l = (n,m)\\ \in \LsetDC \cup \LsetDCRev}} p_{l,t} & \forall n \in \Nset, t \in \Tset, \label{eqn:nodal_balance_p_dc}
\end{align}
with the incidence matrix $K$, the cycle basis $\Cset$, and the cycle incidence matrix $C$ as in \cite{NEUMANN2022118859}.
Eq.\ \labelcref{eqn:kvl} enforces Kirchhoff's voltage law \cite{192975}.
As the active power flows are solely determined by the voltage angle differences and the line reactances, the voltage angle difference limit \labelcref{eqn:voltage_angle_difference} is replaced by \labelcref{eqn:voltage_angle_difference_dc_approx}.
Since this represents a fixed limit on the power flow that does not depend on the number of circuits, we set a tighter upper bound $u_l^{\text{max}}$ for the number of circuits $u_l$ if $u_l^{\text{max}} \geq \theta_l^{\text{max}}/(x_l a_l f_l^{\text{max}})$, which we have found to be necessary in practice to avoid numerical issues.

\subsubsection{DC approximation with transmission losses}\noindent
To account for losses within the DC approximation, a linear relaxation of the loss approximation $p_{l,t}^{\text{loss}} \approx r_l p_{l,t}^2$ can be used \cite{NEUMANN2022118859}.
As before, we ignore all reactive power constraints and replace \labelcref{eqn:ac_flows_p,eqn:ac_nodal_balance_p,eqn:thermal_limit,eqn:voltage_angle_difference,eqn:ac_flows_p} by \labelcref{eqn:kvl}, \labelcref{eqn:voltage_angle_difference_dc_approx} and
\begin{align}
  & p_{l,t}^{\text{loss}} \geq r_l p_{l,0} (2p_{l,t} - p_{l,0}) & \forall l \in \LsetAC, t \in \Tset,p_{l,0} \in \mathcal{P}_{l,0}, \label{eqn:losses_p_dc} \\ 
  & \vert p_{l,t} \vert \leq a_l f_l^{\text{max}} u_l - p_{l,t}^{\text{loss}} & \forall l \in \LsetAC, t \in \Tset, \label{eqn:thermal_limit_dc_lossy} \\
  &\sum_{s \in\Sset_n} p_{s,t} - \sum_{s \in \mathcal{S}_n^{\text{sto}}} p_{s,t}^c - d_{n,t}^p = \sum_{l \in \LsetAC} K_{nl} p_{l,t} + \vert K_{nl} \vert \frac{p_{l,t}^{\text{loss}}}{2} + \sum_{\substack{l = (n,m)\\ \in \LsetDC \cup \LsetDCRev}} p_{l,t} & \forall n \in \Nset, t \in \Tset, \label{eqn:nodal_balance_p_dc_lossy}
\end{align}
with
\begin{equation}\label{eqn:tangent_points}
  \mathcal{P}_{l,0} = \left\{\pm h \cdot \frac{a_l p^{\text{max}}_l u^{\text{max}}_l}{H} \mid h = 1, \dots, H \right\}
\end{equation}
the set of $2H$ points at which the tangents are centred.
Due to the dependence of $\mathcal{P}_{l,0}$ on the maximum number of circuits $u_l^{\text{max}}$, lowering $u_l^{\text{max}}$ to match the voltage angle difference limit as proposed in the previous section changes the feasible space since the loss relaxation is tightened.

\subsubsection{LPAC Approximation}\noindent
The LPAC approximation \cite{lpac,powermodels} assumes that voltage magnitudes are close to their nominal value and voltage angle differences are small.
Furthermore, it relaxes the cosine function within the interval defined by the voltage angle difference limit.
\labelcref{eqn:ac_flows_p,eqn:ac_flows_q} in the original problem are replaced by
\begin{align}
    &p_{l,t} = g_l (v_{n,t} - v_{m,t}) + g_l (1 - \cosrelax_{l,t}) - b_l \theta_{l,t} & \forall l \in \LsetAC \cup \LsetACRev, t \in \Tset, \label{eqn:lpac_p} \\
    &q_{l,t} = -\frac{b_l^{\text{sh}}}{2} (2v_{n,t} - 1) - b_l (v_{n,t} - v_{m,t}) - b_l (1 - \cosrelax_{l,t}) - g_l \theta_{l,t} & \forall l \in \LsetAC \cup \LsetACRev, t \in \Tset, \label{eqn:lpac_q} \\
    &\cosrelax_{l,t} \leq 1 - \frac{1 - \cos(\theta^{\text{max}}_l)}{(\theta^{\text{max}}_l)^2} \theta_{l,t}^2 & \forall l \in \LsetAC \cup \LsetACRev, t \in \Tset, \label{eqn:cosine_relaxation}
\end{align}
where $\cosrelax_{l,t} \in [\cos(\theta^{\text{max}}_l), 1]$ represents the relaxation of the cosine function.
We provide a derivation in \cref{sec:lpac_derivation}.
The formulation implicitly contains losses and reactive power demand of power transmission, for which the respective expressions and their lower bounds are given by
\begin{align}
  p_{l,t} + p_{l^\prime,t} ={}& 2 g_l (1 - \cosrelax_{l,t}) \geq 2 g_l \frac{1 - \cos(\theta_l^{\text{max}})}{(\theta_l^{\text{max}})^2} \theta^2_{l,t}, \label{eqn:lpac_p_losses} \\
  \begin{split}
  q_{l,t} + q_{l^\prime,t} ={}& - b_l^{\text{sh}} (v_{n,t} + v_{m,t} - 1) - 2 b_l (1 - \cosrelax_{l,t}) \\
                            \geq{}& - b_l^{\text{sh}} (v_{n,t} + v_{m,t} - 1) - 2 b_l \frac{1 - \cos(\theta_l^{\text{max}})}{(\theta_l^{\text{max}})^2} \theta^2_{l,t}, \end{split} \label{eqn:lpac_q_losses}
\end{align}
which we derive with \labelcref{eqn:cosine_relaxation} and $g_l \geq 0, b_l \leq 0$.

\subsubsection{Decoupled approximation}\noindent
What we call the \emph{decoupled} approximation represents an adaptation of a family of similar approximations published, e.\,g., in \cite{6495738,fortenbacher_linearquadratic_2019}.
It assumes that voltage magnitudes are close to their nominal value and voltage angles differences are small.
It further approximates transmission losses as $p_{l,t}^{\text{loss}} \approx g_l \theta_{l,t}^2$ and reactive power demand as $q_{l,t}^{\text{dem}} \approx -b_l \theta_{l,t}^2$ and distributes them equally to both endpoints.
We replace \labelcref{eqn:ac_flows_p,eqn:ac_flows_q,eqn:ac_nodal_balance_p,eqn:ac_nodal_balance_q,eqn:thermal_limit} in the original problem by
\begin{align}
    &p_{l,t} = g_l (v_{n,t} - v_{m,t}) - b_l \theta_{l,t} & \forall l \in \LsetAC \cup \LsetACRev, t \in \Tset, \label{eqn:decoupled_p} \\
    &q_{l,t} = - \frac{b_l^{\text{sh}}}{2} (2 v_{n,t} - 1) - b_l (v_{n,t} - v_{m,t}) - g_l \theta_{l,t} & \forall l \in \LsetAC \cup \LsetACRev, t \in \Tset, \label{eqn:decoupled_q}
\end{align}
and the nodal balance equations
\begingroup
\allowdisplaybreaks
\begin{align}
    &\sum_{s \in\Sset_n} p_{s,t} - \sum_{s \in \mathcal{S}_n^{\text{sto}}} p_{s,t}^c - d_{n,t}^p = \sum_{\substack{l \in \LsetAC}} K_{nl} p_{l,t} + \vert K_{nl} \vert \frac{p_{l,t}^{\text{loss}}}{2} + \sum_{\substack{l=(n,m) \\ \in \LsetDC \cup \LsetDCRev}} p_{l,t} & \forall n \in \Nset, t \in \Tset, \label{eqn:ac_nodal_balance_p_decoupled} \\
    &\sum_{s \in\Sset_n} q_{s,t} - d_{n,t}^q = \sum_{\substack{l=(n,m) \\ \in \Lset \cup \Lset}} q_{l,t} + \sum_{l \in \LsetAC} \vert K_{nl} \vert \frac{q_{l,t}^{\text{dem}}}{2} & \forall n \in \Nset, t \in \Tset. \label{eqn:ac_nodal_balance_q_decoupled}
\end{align}
\endgroup
Furthermore, we deviate from \cite{6495738,fortenbacher_linearquadratic_2019}, where the losses are modelled in a mixed-integer formulation, and relax the loss terms via the constraints
\begin{align}
  g_l \theta_{l,t}^2 &\leq p_{l,t}^{\text{loss}} \leq g_l (\theta_l^{\text{max}})^2 & \forall l \in \LsetAC, t \in \Tset, \label{eqn:decoupled_p_losses} \\
  -b_l \theta_{l,t}^2 &\leq q_{l,t}^{\text{dem}} \leq -b_l (\theta_l^{\text{max}})^2 & \forall l \in \LsetAC, t \in \Tset, \label{eqn:decoupled_q_losses}
\end{align}
such that the problem remains convex.
Again deviating from \cite{6495738}, we consider the losses and reactive power demand in the thermal limit constraint 
\begin{align} \label{eqn:thermal_limit_decoupled}
  &(p_{l,t} + p_{l,t}^{\text{loss}})^2 + (q_{l,t} + q_{l,t}^{\text{dem}})^2 \leq (a_l f^{\text{max}}_l u_l)^2 & \forall l \in \LsetAC \cup \LsetACRev, t \in \Tset.
\end{align}

\subsubsection{Comparison of loss approximations}\noindent
Using the relation $r_l / x_l^2 \approx g_l$ for HVAC lines \cite{NEUMANN2022118859}, we find that $p_{l,t}^{\text{loss}} = r_l p_{l,t}^2 = r_l (\theta_{l,t} / x_l)^2 \approx g_l \theta_{l,t}^2$, i.\,e.\ that the losses in \textsc{dc-lossy} and \textsc{decoupled} are approximately equivalent.
The losses in \textsc{lpac} given in \labelcref{eqn:lpac_p_losses} contain an additional term
\begin{equation}
2 \frac{1 - \cos(\theta_l^{\text{max}})}{(\theta_l^{\text{max}})^2} \overset{\theta_l^{\text{max}} \leq \pi/6}{<} 0.9774, \label{eqn:lpac_p_losses_factor}
\end{equation}
which decreases with a larger voltage angle difference limit $\theta_l^{\text{max}}$.
We have bounded the term in \labelcref{eqn:lpac_p_losses_factor} using a typical value for normal operation of  $\theta_l^{\text{max}} = \pi / 6$.
Hence, for the same voltage angle difference, the losses in \textsc{lpac} are more than \qty{2.2}{\percent} lower compared to \textsc{dc-lossy} and \textsc{decoupled}.
Moreover, both \textsc{dc-lossy} and \textsc{decoupled} distribute the losses equally to each endpoint of a branch, while in the \textsc{lpac} approximation the sending endpoint has to compensate them.

\begin{table*}[t]
  \footnotesize
  \caption{Overview of different \textsc{cep} formulations. All formulations use \labelcref{eqn:obj} as the objective function.}\label{tab:model_definitions}
  \centering
    \begin{tabular}{lll}
      \toprule
      problem & type & constraints \\ \midrule
      \textsc{ac-cep} & \makecell[l]{nonconvex nonlinear program} & \labelcref{eqn:online_limit,eqn:online_consistency,eqn:online_boundary,eqn:ps_injections_online,eqn:storage_unit_charging,eqn:storage_unit_complementarity_relaxation,eqn:state_of_charge_limits,eqn:state_of_charge,eqn:state_of_charge_cyclic,eqn:qs_injection_online,eqn:pq_upper,eqn:pq_lower,eqn:ac_flows_p,eqn:ac_flows_q,eqn:thermal_limit,eqn:voltage_angle_difference,eqn:flow_hvdc_along,eqn:flow_hvdc_against,eqn:flow_hvdc_limit,eqn:ac_nodal_balance_p,eqn:ac_nodal_balance_q} \\
      \textsc{dc-cep} & linear program &  \labelcref{eqn:online_limit,eqn:online_consistency,eqn:online_boundary,eqn:ps_injections_online,eqn:storage_unit_charging,eqn:storage_unit_complementarity_relaxation,eqn:state_of_charge_limits,eqn:state_of_charge,eqn:state_of_charge_cyclic,eqn:flow_hvdc_along,eqn:flow_hvdc_against,eqn:flow_hvdc_limit,eqn:kvl,eqn:thermal_limit_dc,eqn:voltage_angle_difference_dc_approx,eqn:nodal_balance_p_dc} \\
      \textsc{dc-lossy-cep} & linear program & \labelcref{eqn:online_limit,eqn:online_consistency,eqn:online_boundary,eqn:ps_injections_online,eqn:storage_unit_charging,eqn:storage_unit_complementarity_relaxation,eqn:state_of_charge_limits,eqn:state_of_charge,eqn:state_of_charge_cyclic,eqn:flow_hvdc_along,eqn:flow_hvdc_against,eqn:flow_hvdc_limit,eqn:losses_p_dc,eqn:kvl,eqn:thermal_limit_dc_lossy,eqn:voltage_angle_difference_dc_approx,eqn:nodal_balance_p_dc_lossy} \\
      \textsc{lpac-cep} & \makecell[l]{convex quadratically constrained program} & \labelcref{eqn:online_limit,eqn:online_consistency,eqn:online_boundary,eqn:ps_injections_online,eqn:storage_unit_charging,eqn:storage_unit_complementarity_relaxation,eqn:state_of_charge_limits,eqn:state_of_charge,eqn:state_of_charge_cyclic,eqn:qs_injection_online,eqn:pq_upper,eqn:pq_lower,eqn:thermal_limit,eqn:voltage_angle_difference,eqn:flow_hvdc_along,eqn:flow_hvdc_against,eqn:flow_hvdc_limit,eqn:ac_nodal_balance_p,eqn:ac_nodal_balance_q,eqn:cosine_relaxation,eqn:lpac_p,eqn:lpac_q} \\
      \textsc{decoupled-cep} & \makecell[l]{convex quadratically constrained program} & \labelcref{eqn:online_limit,eqn:online_consistency,eqn:online_boundary,eqn:ps_injections_online,eqn:storage_unit_charging,eqn:storage_unit_complementarity_relaxation,eqn:state_of_charge_limits,eqn:state_of_charge,eqn:state_of_charge_cyclic,eqn:qs_injection_online,eqn:pq_upper,eqn:pq_lower,eqn:voltage_angle_difference,eqn:flow_hvdc_along,eqn:flow_hvdc_against,eqn:flow_hvdc_limit,eqn:ac_nodal_balance_p_decoupled,eqn:ac_nodal_balance_q_decoupled,eqn:decoupled_p,eqn:decoupled_q,eqn:decoupled_p_losses,eqn:decoupled_q_losses,eqn:thermal_limit_decoupled}\\
      \bottomrule
    \end{tabular}
\end{table*}

\subsubsection{Iterative approach to transmission expansion planning with voltage angle difference limits}\noindent
To remove the nonconvexity introduced by expansion-related change of transmission line parameters, we use an established successive convex programming approach given in \cref{alg:iterative_cep} \cite{HAGSPIEL2014654,PyPSA}.
The approach repeatedly solves the \textsc{cep} and fixes the branch parameters such that they match the transmission expansion of the previous iteration.
Convergence is reached if the difference in transmission expansion to the previous iteration is within a small tolerance.

Fixing branch parameters in conjunction with the voltage angle difference limit \labelcref{eqn:voltage_angle_difference} introduces a possible issue, as it can prevent transmission expansion if the thermal limit constraint is rendered non-binding.
We derive a condition based on the active power flow in the DC approximation, which we take as an upper bound for the active power flow in the power flow approximations included in this publication.
While this is not an analytical bound, the DC approximation generally overestimates the feasible power flow from a system perspective as losses are ignored.
With \labelcref{eqn:thermal_limit_dc,eqn:voltage_angle_difference_dc_approx} we obtain
\begin{equation} \label{eqn:condition_vang_diff}
  \vert p_{l,t} \vert \leq a_l f_l^{\text{max}} u_l^k \leq a_l f_l^{\text{max}} u_l^{\text{min}} \Leftrightarrow u_l^{\text{min}} \geq \frac{\theta_l^{\text{max}}}{x_l^k a_l f_l^{\text{max}}},
\end{equation}
where $u_l^k$, $x_l^k$ represent the number of circuits and reactance of line $l \in \LsetAC$ in iteration $k$ of \cref{alg:iterative_cep}, respectively.
If this condition is fulfilled, the voltage angle difference limit prevents transmission expansion.
This can occur, e.\,g., due to the use of open data sources, where series compensation of transmission lines is often ignored \cite{HORSCH2018207}.
We evaluate the condition in the numerical experiment to assess whether the problem instances are affected by this issue.

\subsection{AC-reinforcement} \label{sec:ac_reinforcement}\noindent
We propose a cost-minimizing AC-reinforcement approach similar to \cite{6308747}, which reinforces an initial solution such that it becomes AC-feasible.
The approach presented in \cref{alg:ac_reinforcement} determines a solution of an \textsc{ac generation expansion problem} (\textsc{ac-gep}), where existing power sources can be redispatched and additional capacities can be installed, e.\,g.\ to provide voltage support.
To keep the \textsc{ac-gep} instance tractable, we decompose it into its individual snapshots, which we then solve separately.
Based on the initial solution, we fix transmission and storage system expansion and bound the expansion of the remaining power sources from below.
We use setpoints from the initial solution for $\beta_{s,t}$ and $\beta_{s,t}^{\text{su}}$ such that deviations from the schedule are only possible at an additional cost.
To capture shutdown procedures, we introduce the variable $\beta_{s,t}^{\text{sd}}$.
A penalty $o_s^ {\text{su}} \beta_{s,t}^{\text{sd}}$ is added in the objective function to account for the additional startup costs required to match the subsequent setpoint $\beta_{s,t+1}$.
All variables regarding storage systems are fixed except for $p_{s,t}$, for which we set an upper bound based on the initial solution as in \cite{NEUMANN2022118859}.

In \cref{alg:ac_reinforcement}, we first determine all snapshots for which a reinforcement is necessary by attempting to find a solution to the \textsc{ac-opf} problem obtained after fixing $u_s = u_s^{\star}$ for all $s \in \Sset$.
Thereby, we aim to avoid convergence of a heuristic solver to local minima involving reinforcements that are not necessary.
Furthermore, since the remaining inter-temporal dependencies are resolved in \textsc{ac-opf}, this step can be parallelised.
All snapshots for which no feasible solution could be found are gathered in the set $\mathcal{T}^\prime$.
We then reinforce these snapshots in a sequential procedure where we iterate over $\mathcal{T}^\prime$ in chronological order.
In each iteration, we solve the \textsc{ac-gep} for a single snapshot and adopt the reinforcements for the remaining snapshots.
Here, we again first attempt to find a solution for \textsc{ac-opf}.
\begin{algorithm}[t]
  \caption{AC-reinforcement}\label{alg:ac_reinforcement}
  \begin{algorithmic}
    \STATE \textbf{Input:} Initial solution $\mathfrak{X}^\star$ of \textsc{cep}.  
    \STATE \textbf{Output:} AC-feasible solution $\mathfrak{X}^\prime$ of \textsc{cep}.
    \STATE Define the \textsc{ac-gep}
      \begin{align*}
        \min       \quad  & \mathrlap{\sum_{s \in\Sset} \left( c_s u_s + \sum_{t \in \Tset} \delta_{t} o_s p_{s,t} + o_s^{\text{su}} \left( \beta_{s,t}^{\text{su}} + \beta_{s,t}^{\text{sd}} \right)  \right)} \\
        \text{s.t.} \quad & \mathrlap{\text{\labelcref{eqn:online_limit,eqn:ps_injections_online,eqn:qs_injection_online,,eqn:pq_lower,eqn:pq_upper,eqn:ac_flows_p,eqn:ac_flows_q,eqn:thermal_limit,eqn:voltage_angle_difference,eqn:flow_hvdc_against,eqn:flow_hvdc_along,eqn:flow_hvdc_limit,eqn:ac_nodal_balance_p,eqn:ac_nodal_balance_q}}} & \\
                      & u_l = u_l^{\star} & \forall l \in \Lset \\
                      & u_s = u_s^{\star} & \forall s \in \SsetSto \\
                      & u_s \geq u_s^{\star} & \forall s \in \Sset \setminus \SsetSto, \\
                      & \beta_{s,t} = \beta_{s,t-1}^\star + \beta_{s,t}^{\text{su}} - \beta_{s,t}^{\text{sd}}  & \forall s \in \Sset \setminus \SsetSto, t \in \Tset \\
                      & \beta_{s,t}^{\text{su}} \geq \beta_{s,t}^{\text{su},\star} & \forall s \in \Sset \setminus \SsetSto, t \in \Tset \\
                      & \beta_{s,t}^{\text{su}} = \beta_{s,t}^{\text{su},\star}, \; p_{s,t}^c = p_{s,t}^{c,\star} & \forall s \in \SsetSto, t \in \Tset
      \end{align*}
    where all inter-temporal dependencies except for the expansion decisions are resolved based on assigned values in $\mathfrak{X}^\star$.
    Let \textsc{ac-opf}$_t$ be the \textsc{ac-opf} problem resulting from fixing $u_s = u_s^\star$ for all $s \in \Sset$ and $\Tset = \{t\}$ in the \textsc{ac-gep}.
    \STATE $\mathcal{T}^\prime \leftarrow \left\{\; t \mid \text{no solution found for } \textsc{ac-opf}_t \right\}$
    \STATE \textbf{for} $t \in \mathcal{T}^\prime$ \textbf{do}'
    \STATE \hspace{0.5cm} Attempt to find a solution for \textsc{ac-opf}$_t$.
    \STATE \hspace{0.5cm} \textbf{if} no solution could be found for \textsc{ac-opf}$_t$
    \STATE \hspace{1cm} Determine a solution for \textsc{ac-gep} with $\Tset = \{t\}$.
    \STATE \hspace{1cm} Update $u_s^\star$ based on the reinforcements.
    \STATE \hspace{0.5cm} \textbf{end if}
    \STATE \textbf{end for}
  \end{algorithmic}
\end{algorithm}
The AC-feasibility of the resulting system directly follows from the AC-feasibility of each individual snapshot.

\section{Numerical Experiment}\noindent
We define the two scenarios \emph{SG-dominated} and \emph{IBR-dominated} which represent systems with a high share of SGs and a high share of IBRs, respectively.
For each scenario and each of the power flow approximations \textsc{dc}, \textsc{dc-lossy}, \textsc{lpac}, \textsc{decoupled}, we solve the corresponding \textsc{cep} via \cref{alg:iterative_cep} and reinforce the solution for AC-feasibility using \cref{alg:ac_reinforcement}.
We do not consider decommissioning.

As a data source, we use PyPSA-Eur (v2025.07.0) \cite{HORSCH2018207}, from which we extract the islands of Great Britain and Ireland at their full resolution of \num{354} buses and \num{478} transmission lines.
This allows us to exclude effects of spatial aggregation, which would be required for problem instances of larger spatial extent due to computational constraints.
Furthermore, we aggregate the data temporally to 4380 segments based on \cite{7982024} to remain within reasonable computational limits.
We adopt all assumptions for capital and operational costs from \cite{HORSCH2018207}, unless explicitly mentioned below.
All transmission lines are assumed to be of the type \qty{380}{\kilo\volt} 4-bundle 240/40 Al/St \cite{oeding_elektrische_2011}.
Transformers are not considered.
We adopt the startup costs for SGs from \cite{11305559}.
As in \cite{NEUMANN2022118859}, our base assumption is that the number of circuits of HVAC branches can be at most doubled or increased by three.
We set a voltage angle difference limit of $\theta_l^{\text{max}} = \pi / 6$ to remain within normal operation limits and choose a commonly used conservative line utilization limit of $a_l=0.7$ for all HVAC lines.
The voltage angle difference limit imposes an implicit limit on the transmission expansion, which might be more restrictive than the thermal limit when ignoring the expansion-related change of the line impedance in each iteration of \cref{alg:iterative_cep}.
To consider if this prevents transmission expansion, we evaluate \labelcref{eqn:condition_vang_diff} for all $l \in \LsetAC$ and find that the condition is not fulfilled for any $l$.
All HVAC transmission lines in the grid topology of our problem instances can be extended.
For all except one, at least a doubling of the number of circuits is possible.
HVDC branches can be increased by up to an equivalent of \qty{20}{\giga\watt}.
We set the loss factor $\eta_l$ for HVDC transmission to \qty{3}{\percent} per \qty{1000}{\kilo\meter} and assume that new HVDC circuits are connected via voltage-source converter stations.
To define the loss relaxation in \textsc{dc-lossy-cep}, we set $H=3$ in \labelcref{eqn:tangent_points} as in \cite{NEUMANN2022118859}.

In both scenarios, we allow expansion of solar and wind power and flexibly controllable voltage-increasing and voltage-decreasing reactive power compensation devices.
The capital costs of \qty{20}{\euro\per{\kilo\volt\ampere}} for voltage-increasing and \qty{26}{\euro\per{\kilo\volt\ampere}} for voltage-decreasing compensation devices correspond to those of \emph{mechanically-switched capacitors with damping network} and \emph{shunt reactors} in the German transmission expansion plan \cite{nep_kosten}.
We assume a lifetime of 20 years for both.
As we ignore switching due to the unit commitment relaxation, we do not include further types of compensation devices whose advantages would not be considered in our model.
In \sgdominated, we additionally allow expansion of combined-cycle gas turbines (CCGT).
In \ibrdominated, all existing conventional generation except hydropower is removed as in\cite{NEUMANN2022118859} and, additionally, expansion of battery energy storage systems (BESS) is allowed.

We furthermore provide a set of assumptions for modelling reactive power. 
Reactive loads are given by assuming an inductive grid with a fixed power factor of $\cos \varphi_d = 0.99$ as in \cite{7232601}.
To model the reactive power capabilities of all power sources, we provide four possible conservative approximations of PQ-capability curves as given in \cref{tab:pq_curves}.
Synchronous generators are assigned the \emph{d-curve} capability curve, offshore wind inverters are assigned the \emph{u-shape} curve, voltage-source HVDC converters and BESS are assigned the \emph{rectangle} variant.
All other power sources are assumed to have the \emph{triangle} capability curve, which models a typical requirement from grid codes for lower voltage levels of requiring a power factor range of $\cos \varphi_s = 0.95$ both inductive and capacitive.  
Through the definition of the PQ capabilities, we ensure that the apparent load of an IBR is within the typical rating of inverters for a power of $1.1 \, p_s^{\text{max}}$ \cite{uk_qflex}.
The apparent power limit $\sqrt{p_{s,t}^2 + q_{s,t}^2} \leq 1.1 \, p_s^{\text{max}}$ holds as $q_{s,t} \leq 0.46 \, \text{(p.u.)}$ in all PQ capability curve approximations.
Moreover, \labelcref{eqn:storage_unit_complementarity_relaxation} ensures that simultaneous charging and discharging of a storage system remains below $p_s^{\text{max}} u_s$, hence no special treatment of BESS inverters is necessary.
We do not include other inverter options, which may be able to provide more reactive power.
As an equivalent, it is possible to install voltage-increasing and voltage-decreasing compensation at the respective bus.
\begin{table*}[t]
  \footnotesize
  \caption{PQ-capability curve approximations used in this publication \cite{uk_qflex}. For \emph{triangle} we assume $\cos \varphi = 0.95$.}\label{tab:pq_curves}
  \centering
  \begin{tabular}{lllll}
    \toprule
    & $\mathcal{K}^{\text{up}}_s$ & $\mathcal{K}^{\text{lo}}_s$ & $q^{\text{min}}_s$ & $q^{\text{max}}_s$   \\ \midrule
    \emph{d-curve} & $\{(\frac{1}{2}, 1), (-\frac{1}{3}, 1)\}$ & $\emptyset$ & $-0.4$ & $0.6$ \\
    \emph{u-shape} & $\emptyset$ & $\{(\frac{1}{2}, 0), (-\frac{1}{2}, 0)\}$ & $-0.4$ & $0.4$ \\
    \emph{triangle} & $\emptyset$ & $\{(\frac{1}{\tan \varphi}, 0), (-\frac{1}{\tan \varphi}, 0) \}$ & $-\tan \varphi$ & $\tan \varphi$ \\
    \emph{rectangle} & $\emptyset$ & $\emptyset$ & $-0.4$ & $0.4$ \\ \bottomrule
  \end{tabular}
\end{table*}

All \textsc{cep} instances are defined using GAMSPy and solved using Gurobi 12.0.3 on a compute node of the German Aerospace Center's supercomputer CARO, which is equipped with an AMD EPYC 7702 processor and one Terabyte of main memory.
All \textsc{ac-opf} and \textsc{ac-gep} problem instances are solved with IPOPT \cite{wachter_implementation_2006}.

\section{Results}\noindent
\begin{figure}[!t]
\centering
    \includegraphics{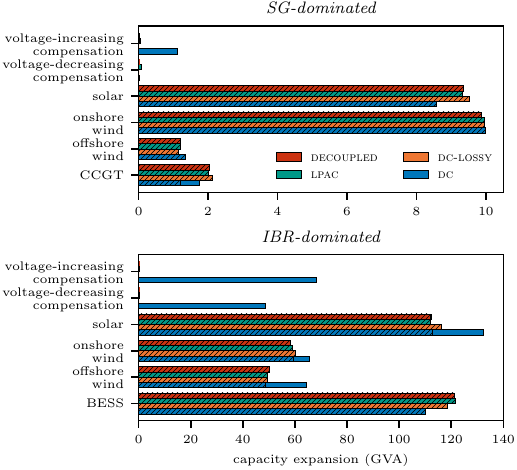}
    \caption{Capacity expansion decisions for each scenario and power flow approximation, excluding transmission expansion. Hatched areas indicate the respective value in the initial solution before the AC-reinforcement.} \label{fig:capacities}
\end{figure}
\begin{figure}[!t]
\centering
    \includegraphics{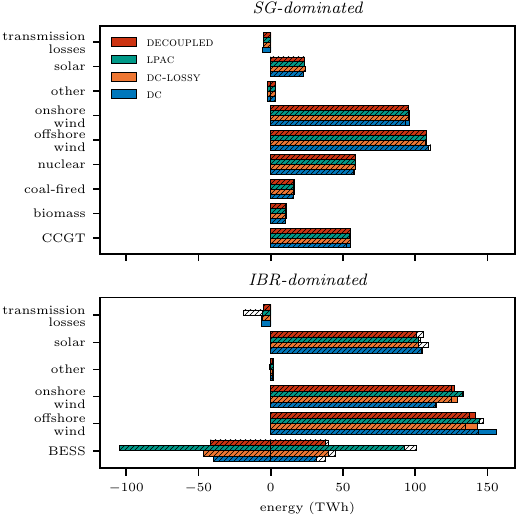}
    \caption{Energy mix for each scenario and power flow approximation. Hatched areas indicate the respective value in the initial solution before the AC-reinforcement. The total load is equal to \qty{362}{\tera\watt\hour}.} \label{fig:energy}
\end{figure}
\begin{figure*}[!t]
\centering
    \includegraphics[width=\textwidth]{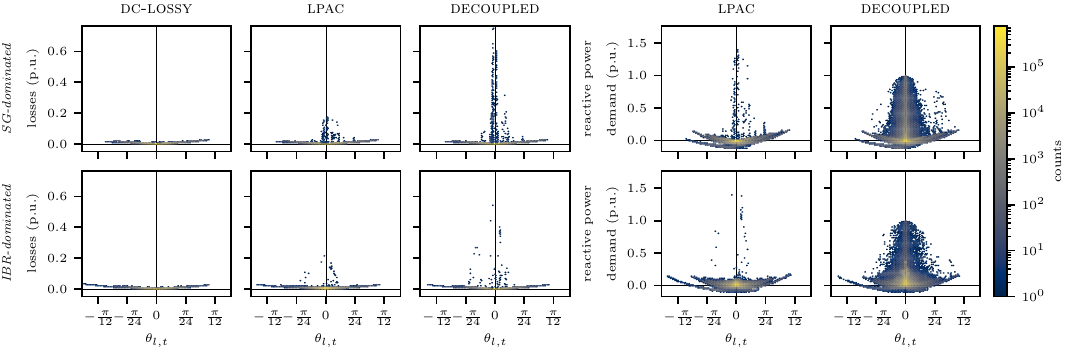}
    \caption{Losses and reactive power demand of transmission in p.u.\ of the thermal limit as a function of the voltage angle difference. Losses and reactive power demand are given by the sum of the active and reactive power flows at both endpoints of each branch, respectively. We obtain the voltage angle differences for \dclossy{} via the relation $\theta_{l,t} = p_{l,t} x_l$.} \label{fig:loss_relaxations}
\end{figure*}
\begin{figure}[!t]
  \centering
  \includegraphics{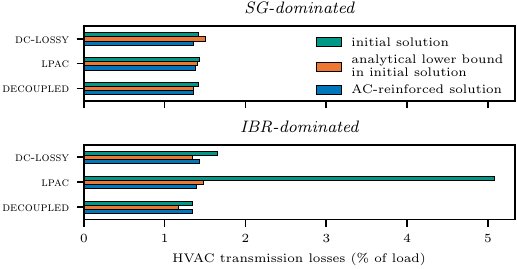}
  \caption{Transmission losses in initial solution in comparison to their analytical lower bound and the losses in the AC-reinforced solution. The lower bounds are computed by evaluating the analytical expressions with the assigned values of the respective decision variables. In the case of \dclossy, the linear relaxation of the loss approximation can lead to losses that are lower than the analytical lower bound.}
    \label{fig:losses_barplot}
\end{figure}

\cref{tab:results} shows summarizing quantities, while \cref{fig:capacities,fig:energy} show the capacity expansion and energy mix for each solution, respectively.
\cref{fig:loss_relaxations,fig:losses_barplot} show the losses for each solution, if modelled.
In the following, we describe the results qualitatively by investigating the capacity mixes of the resulting systems, both before and after the AC-reinforcement.
We reason about the AC-feasibility of the solutions before the reinforcement by investigating the number of \emph{non-converging snapshots}, i.\,e.\ snapshots for which no solution to the \textsc{ac-opf} problem could be found\footnote{While non-convergence of a heuristic solver on an \textsc{ac-opf} problem instance is not a proof of infeasibility, we take it as an indicator for the ill-conditionedness of the problem instance as done e.\,g.\ in \cite{NEUMANN2022118859}.
} in \cref{alg:ac_reinforcement}.
Moreover, we investigate the common issue with power flow approximations of \emph{fictitious}, i.\,e.\ non-physical losses and reactive power consumption of transmission by deviation from the lower bound of the respective relaxation \cite{NEUMANN2022118859}.  
This can occur if increasing losses or reactive power consumption either does not significantly change or decreases the objective value.

All approximations required redispatch or reinforcements in \cref{alg:ac_reinforcement}.
The approach led to fewer corrections for initial solutions based on \dclossy, \lpac, or \decoupled, which all include transmission losses.
For the lossless \dc{} approximation the approach led to a distorted solution in the \ibrdominated{} scenario due to the high share of non-converging snapshots.
In \sgdominated, \dc{} shows a significantly larger amount of redispatch.
For both scenarios, this leads to the highest system costs.
We hence exclude these results in the following.

The differences in system costs between the remaining approximations are negligible, both before and after the reinforcement.
Some differences exist in the expansion decisions, the share of AC-feasible snapshots before the reinforcement, and the required redispatch.
In \sgdominated, there are no significant differences in the expansion decisions as decommissioning was not considered and due to the high share of existing synchronous generation.
In \ibrdominated{} however, \dclossy{} shows at least \qty{5.2}{\giga\watt} higher investments into solar and wind power and \qty{2.9}{\giga\watt} lower investments into BESS compared to \lpac{} and \decoupled, while transmission expansion is around \qty{0.3}{\tera\watt\kilo\meter} higher.
We do not observe improvements over \dclossy{} in terms of AC-feasibility before the AC-reinforcement in \lpac{} or \decoupled.
In \sgdominated, the problem instance is likely almost AC-feasible before the expansion.
In \ibrdominated{}, \dclossy{} shows the best results with close to \qty{100}{\percent} AC-feasible snapshots, but \lpac{} and \decoupled{} perform similarly well.
The reinforcements due to \cref{alg:ac_reinforcement} are negligible in both scenarios.

\cref{fig:loss_relaxations,fig:losses_barplot} show that the relaxations led to fictitious transmission losses and reactive power consumption.
It is visible in \cref{fig:loss_relaxations} that this only occurs for small voltage angle differences, i.\,e.\ if a line is under low load.
As expansion is associated with much higher costs than power injection, the expansion decisions are not affected by fictitious losses.
Furthermore, \cref{fig:losses_barplot} shows that the analytical lower bounds for the transmission losses in all approximations are close to the losses in the AC-reinforced solution.
In the case of \lpac{} and \decoupled{}, fictitious reactive power consumption allows the solver to use the branch as a voltage-decreasing compensation device.
This effect also occurs mostly for small voltage angle differences, although the values are more widely dispersed in comparison to the losses.
Due to the coupling between active and reactive power in \lpac{} via the cosine relaxation term $\cosrelax_{l,t}$, fictitious reactive power consumption occurs if and only if fictitious losses occur.
Since losses are associated with costs, it is reasonable to expect that fictitious reactive power consumption would be lower in \lpac.
However, while it is visible in \cref{fig:loss_relaxations} that the number of severe outliers is reduced in comparison to \decoupled, \cref{fig:losses_barplot} shows significantly higher fictitious losses for \lpac{} in the \ibrdominated{} scenario, which implies higher fictitious reactive power consumption.
We attribute this to the flat optimum typical for \textsc{cep} instances of IBR-dominated power systems \cite{PEDERSEN2021121294}, which in this case leads to small differences in the objective value between installing compensation devices and accepting higher losses to allow fictitious reactive power consumption from branches. 
Decreasing the convergence tolerance of the solver from $10^{-5}$ to $10^{-6}$ already lessened the severity of the issue significantly.
This suggests that relaxed coupling terms between active and reactive power, which are commonly used in power flow approximations \cite{lpac,qc_relaxation,7763860,10646561}, lead to a numerically more challenging problem in the case of a flat optimum.  
Similarly, the significantly higher amount of storage cycling of \lpac{} in \ibrdominated{} is likely also related to the flat optimum.

\cref{tab:results} shows that the computational tractability of \lpac{} and \decoupled{} is significantly worse in comparison to \dclossy, as solution time can increase by up to \num{14} times while maximum memory usage can quadruple.

\begingroup
\setlength{\tabcolsep}{2.95pt}
\begin{table*}[t]
  \footnotesize
\centering
\caption{Summarizing data on system costs and computational tractability of the different formulations.} \label{tab:results}
\begin{tabularx}{\textwidth}{Xrrrrrrrr}
\toprule
 & \multicolumn{4}{l}{\emph{SG-dominated}} & \multicolumn{4}{l}{\emph{IBR-dominated}} \\
 & \textsc{dc} & \textsc{dc-lossy} & \textsc{lpac} & \textsc{decoupled} & \textsc{dc} & \textsc{dc-lossy} & \textsc{lpac} & \textsc{decoupled} \\
\midrule
AC-feasible snapshots of initial solution (\unit{\percent}) & 99.93 & 100.00 & 99.98 & 100.00 & 89.95 & 99.91 & 99.18 & 99.36 \\
\midrule
positive redispatch in \cref{alg:ac_reinforcement} (\unit{\tera\watt\hour}) & 15.51 & 4.88 & 4.27 & 4.56 & 84.01 & 30.86 & 19.40 & 22.46 \\
negative redispatch in \cref{alg:ac_reinforcement} (\unit{\tera\watt\hour}) & -9.41 & -5.11 & -4.45 & -4.76 & -77.75 & -31.41 & -32.67 & -22.30 \\
\midrule
system costs, initial solution ($10^9$ \unit{\euro\per\year}) & 29.28 & 29.48 & 29.47 & 29.47 & 26.82 & 27.54 & 27.52 & 27.53 \\
system costs ($10^9$ \unit{\euro\per\year}) & 29.53 & 29.47 & 29.47 & 29.47 & 31.17 & 27.55 & 27.53 & 27.55 \\
\midrule
transmission expansion (\unit{\tera\watt\kilo\meter}) & 3.14 & 2.34 & 2.45 & 2.41 & 12.29 & 12.05 & 11.60 & 11.72 \\
\midrule
wall-clock time, \cref{alg:iterative_cep} (\unit{\hour}) & 1.75 & 2.20 & 26.90 & 31.70 & 5.60 & 11.40 & 60.20 & 70.17 \\
times solved, \cref{alg:iterative_cep} & 3 & 3 & 4 & 3 & 4 & 4 & 5 & 5 \\
max.\ memory usage, \cref{alg:iterative_cep} (\unit{\giga\byte}) & 67.22 & 104.10 & 193.50 & 416.17 & 93.07 & 87.05 & 276.60 & 380.31 \\
max.\ number of variables, \cref{alg:iterative_cep} ($10^6$) & 5.34 & 8.10 & 36.93 & 51.46 & 9.23 & 11.75 & 41.21 & 55.59 \\
max.\ number of constraints, \cref{alg:iterative_cep} ($10^6$) & 9.93 & 22.85 & 50.76 & 65.30 & 16.42 & 29.16 & 57.39 & 71.97 \\
 \bottomrule
\end{tabularx}
\end{table*}
\endgroup

\section{Discussion}\noindent
Due to non-converging snapshots and necessary redispatch in all initial solutions, a method like our AC-reinforcement approach must be used in any case to ensure AC-feasibility.
For improving the AC-feasibility of the initial solution given problem instances at the level of detail we investigated, considering transmission losses is most decisive.
We do not recommend using the lossless \dc{} approximation for this reason.
We instead recommend \dclossy{} by default, as it offers clear improvements over \dc{} while remaining computationally tractable.
Since \lpac{} and \decoupled{} are less computationally tractable but led to no increase in the number of AC-feasible snapshots, we do not recommend their usage for the type of \textsc{cep} instances considered here.

A shortcoming of our study is that we did not investigate more recent advancements in power flow approximations, which e.\,g.\ provide tighter relaxations \cite{10646561,7763860}.
To model the ability of SGs to provide reactive power compensation, the unit commitment relaxation from \cite{11305559} was necessary.
While it is an improvement over not considering startup procedures at all, as it is common in related studies \cite{NEUMANN2022118859}, the approximation error compared to a true unit commitment solution has not been investigated.
Moreover, we ignore transmission line switching, thereby possibly overestimating the need for voltage-decreasing compensation during low load.  
While this issue is less severe with \lpac{} or \decoupled{}, which are able to offset the capacitive load from the grid  via  fictitious reactive power consumption from transmission lines, this is not the the case for the reinforced \textsc{ac} solution.
Due to continuous expansion decisions and cost minimization, the resulting systems are designed to operate at their limits.
With discrete expansion decisions, some forced overcapacities could possibly result in an improved AC-feasibility of the system.

Our results agree with previous research, which finds that approximate solutions are mostly not AC-feasible \cite{VENZKE2020106480,6308747}.  
Moreover, previous research also finds improvements regarding AC-feasibility or the amount of necessary reinforcements when using approximations that include transmission losses over the \dc{} approximation \cite{NEUMANN2022118859,6308747,WOGRIN2020115925,DAVE2021106683,7038299}.
In contrast to \cite{NEUMANN2022118859}, some non-converging snapshots remain in our experiment.
This can be attributed to the higher spatial resolution considered, which exposes voltage problems, and the weaker assumptions used.
The authors of \cite{WOGRIN2020115925} find that reactive power constraints are imperative for improving AC-feasibility.
Our results suggest that it is rather the inclusion of transmission losses that brings the largest improvement.
However, as the authors did not investigate \dclossy{} and due to further differences in the scope of the model and the experiment, we do not consider this a contradiction.

\section{Conclusion and Outlook}\noindent
We formulated a multi-period alternating current capacity expansion problem that considers operational constraints and simplified stability limits.
We obtain solutions to this problem by redispatching and reinforcing an initial solution determined by solving a convex approximation of the original problem.
Our findings show that our approach requires less redispatch and reinforcements if the initial solution is determined using a power flow approximation that includes transmission losses.
We also find transmission losses to be the most decisive factor for improving the AC-feasibility of the initial solution, although redispatch or reinforcements were necessary for each of the analysed power flow approximations.
Considering reactive power constraints in the problem formulation for the initial solution did not lead to further improvements, but significantly reduced the computational tractability.
Due to non-converging snapshots, fictitious losses and fictitious reactive power consumption in the initial solutions, we recommend to always use an a posteriori AC-reinforcement approach.

Future work can use the resulting operating points as initial conditions for dynamic grid simulations.
Furthermore, the quality of the reinforced solution could be assessed by using a relaxation of the power flow equations such as \cite{qc_relaxation}, which would allow obtaining a lower bound to the costs of an AC-feasible system.
We have moreover identified an issue that arises when using voltage angle difference constraints within the iterative transmission expansion planning, especially for strongly aggregated networks that are common in studies considering areas of large spatial extent.
Finally, we observed that a power flow approximation with a relaxed coupling term between active and reactive power can be numerically more challenging for a problem instance with a flat optimum, leading to increased occurrence of undesired effects such as fictitious losses.

\section*{Acknowledgments}
The authors gratefully acknowledge funding from the German Federal Ministry for Economic Affairs and Energy under grant numbers FKZ 03EI1055A, FKZ 03EI1055B and the scientific support and high performance computing (HPC) resources provided by the German Aerospace Center (DLR).
The HPC system CARO is partially funded by the ``Ministry of Science and Culture of Lower Saxony'' and the ``Federal Ministry for Economic Affairs and Energy''.
GR thanks M.\ Wetzel (DLR) for repeated exchanges.

\section*{Appendix}
\renewcommand{\thealgorithm}{A.\arabic{algorithm}}
\setcounter{algorithm}{0}

\label{sec:lpac_derivation} \noindent
The \lpac{} approximation was originally published in \cite{lpac} and later improved in \cite{powermodels}.
We provide an updated derivation here.
Assuming $v_n, v_m \approx 1$, we obtain from \labelcref{eqn:ac_flows_p,eqn:ac_flows_q}
\begin{align*}
  p_l' & = g_l -  g_l \cos\theta_l - b_l \sin \theta_l, \\
  q_l' & = - (b_l + \frac{b^\text{sh}_l}{2}) - g_l \sin \theta_l + b_l \cos\theta_l.
\end{align*}
We consider the error introduced by this step via the linear terms of a Taylor expansion of $p_l^{\text{ac}} - p_l'$, $q_l^{\text{ac}} - q_l'$ around the flat voltage profile $v_{n,0},v_{m,0}=1$, $\theta_{l,0}=0$ given by
\begin{align*}
  p_l^\Delta = p_l^{\text{ac}} - p_l' & \approx  2 g_l (v_n - 1) - g_l(v_n + v_m - 2),                              \\
  q_l^\Delta = q_l^{\text{ac}} - q_l' & \approx - 2 (b_l - \frac{b^\text{sh}_l}{2}) (v_n - 1) - b_l(v_n + v_m - 2),
\end{align*}
such that $p_l^{\text{ac}} \approx p_l' + p_l^\Delta$, $q_l^{\text{ac}} \approx q_l' + q_l^\Delta$.
With $\sin \theta_l \approx \theta_l$ and $\cosrelax_l$ as in \labelcref{eqn:cosine_relaxation}, we obtain \labelcref{eqn:lpac_p,eqn:lpac_q}.

\begin{algorithm}[H]
  \caption{Capacity expansion planning via successive convex programming \cite{HAGSPIEL2014654,PyPSA}}
  \label{alg:iterative_cep}
  \begin{algorithmic}
    \STATE $u_l^{\star,0} \leftarrow u^{\text{min}}_l$   
    \STATE $k = 0$  
    \STATE $\Delta = \infty$  
    \STATE \textbf{while} $\Delta > 0.05$ \textbf{do}
    \STATE \hspace{0.5cm} $k \leftarrow k + 1$
    \STATE \hspace{0.5cm} $u^{\star,k}_l \leftarrow$ Solve convex formulation of capacity expansion problem with parameters of HVAC branches $\LsetAC$ fixed to $\underline{y}_l^k = \underline{y}_l^{\prime} u_{l}^{\star,k-1}$, $b_l^{\text{sh},k} = b_l^{\text{sh}\prime} u_{l}^{\star,k-1}$, where $\underline{y}_l^{\prime}$, $b_l^{\text{sh}\prime}$ are given per circuit.
    \STATE \hspace{0.5cm} $\Delta \leftarrow \frac{\lVert u_l^{\star,k} - u_l^{k-1,\star} \rVert_2}{\lVert u_l^{\star,k} \rVert_2}$
    \STATE \textbf{end while}
    \STATE Fix the number of circuits of HVAC branches $\LsetAC$ to $u_l^{\star,k}$, determine the corresponding branch parameters and solve the capacity expansion problem.
  \end{algorithmic}
\end{algorithm}

\bibliography{bibliography}
\bibliographystyle{abbrv}

\end{document}

%% file: shortcuts.tex
\newcommand{\Sset}{\mathcal{S}}
\newcommand{\SsetSto}{\mathcal{S}^{\text{sto}}}
\newcommand{\SsetSg}{\mathcal{S}^{\text{sg}}}
\newcommand{\SsetIbr}{\mathcal{S}^{\text{ibr}}}
\newcommand{\Tset}{\mathcal{T}}
\newcommand{\Nset}{\mathcal{N}}
\newcommand{\Lset}{\mathcal{L}}
\newcommand{\LsetRev}{\mathcal{L}^\prime}
\newcommand{\LsetAC}{\mathcal{L}^{\text{ac}}}
\newcommand{\LsetACRev}{\mathcal{L}^{\text{ac}\prime}}
\newcommand{\LsetDC}{\mathcal{L}^{\text{dc}}}
\newcommand{\LsetDCRev}{\mathcal{L}^{\text{dc}\prime}}

\newcommand{\Cset}{\mathcal{C}}

\newcommand{\jim}{\mathrm{j}} 

\DeclareMathOperator{\cosrelax}{\breve{\cos}}

\newcommand{\dc}{\textsc{dc}}
\newcommand{\dclossy}{\textsc{dc-lossy}}
\newcommand{\lpac}{\textsc{lpac}}
\newcommand{\decoupled}{\textsc{decoupled}}

\newcommand{\sgdominated}{\emph{SG-dominated}}
\newcommand{\ibrdominated}{\emph{IBR-dominated}}